\documentclass[a4paper,leqno]{article}

\usepackage{amssymb}
\date{\empty}

\newcommand{\ord}{\,{\rm ord}\,}
\newcommand{\ccc}{{\mathbb{C}}}
\newcommand{\rr}{{\mathbb{R}}}
\newcommand{\nn}{{\mathbb{N}}}

\newtheorem{Twierdzenie}{Theorem}[section]

\newtheorem{Lemat}[Twierdzenie]{Lemma}

\newenvironment{Example}{\vspace{1em}\noindent{\bf Example\refstepcounter{Twierdzenie}\ \theTwierdzenie\ }}{\vspace{1em}}

\newenvironment{Uwaga}{\vspace{1em}\noindent{\bf Remark
\refstepcounter{Twierdzenie}\theTwierdzenie\ }}{\vspace{1em}}

\begin{document}
\newtheorem{Wlasnosc}[Twierdzenie]{Property}
\title{A bound for the Milnor number of plane curve singularities}
\author{Arkadiusz P\l oski}
\maketitle

\noindent {\bf Abstract}. Let $f=0$ be a plane algebraic curve of degree $d>1$ with an isolated singular point at $0\in\ccc^2$. We show that the Milnor number $\mu_0(f)$ is less than or equal to $(d-1)^2-\left[\frac{d}{2}\right]$, unless $f=0$ is a set of $d$ concurrent lines passing through $0$. Then we characterize the curves $f=0$ for which $\mu_0(f)=(d-1)^2-\left[\frac{d}{2}\right]$.\footnotetext{2010 Mathematics Subject Classification: 14B05, 14N99.\\
\hspace*{4.5ex}Keywords: Milnor number, plane algebraic curve.}

\addtocounter{section}{-1}
\section{Introduction}
Let $f\in\ccc[x,y]$ be a polynomial of degree $d>1$ such that the curve $f=0$ has an isolated singularity at the origin $0\in\ccc^2$. Let $\mathcal{O}=\mathcal{O}_{\ccc^2,0}$ be the ring of germs of holomorphic functions at $0\in\ccc^2$. The Milnor number $\mu_0(f)=\dim_{\ccc}{}^{\mathcal{O}}\!\!\!\diagup_{\!\!\left(\frac{\partial f}{\partial x},\frac{\partial f}{\partial y}\right)}$ is less than or equal to $(d-1)^2$ by B\' ezout's theorem. The equality $\mu_0(f)=(d-1)^2$ holds if and only if $f$ is a homogeneous polynomial. The aim of this note is to determine the maximum Milnor number $\mu_0(f)$ for non-homogeneous polynomials $f$ (Theorem \ref{thm1.1}) and to characterize the polynomials for which this maximum is attained (Theorem \ref{thm1.4}). The general problem to describe singularities that can occur on plane curves of given degree was studied by G.~M.~Greuel, C.~Lossen and E.~Shustin in \cite{Greuel} (see also \cite{Wall}, Chapter 7 for further references). A bound for the sum of the Milnor numbers of projective hypersurfaces with isolated singular points was given recently by June Huh in \cite{Huh}. Note here that a result of this type follows from Pl\" ucker-Teissier's formula for the degree of the dual hypersurface (see \cite{Teissier}, Appendix 2).

Let us recall usual notions and conventions. By the curve $f=0$ we mean (see \cite{Fulton}, Chapter 3) the linear subspace $\ccc\, f$ of $\ccc[x,y]$. If the polynomial $f$ has no multiple factors then we identify the curve $f=0$ and the set $\{P\in\ccc^2:f(P)=0\}$. We denote by $\ord_0 f$ the order of the polynomial $f$ at $0\in\ccc^2$ and by $i_0(f,g)=\dim_\ccc{}^{\mathcal{O}}\!\!\!\diagup_{\!\!(f,g)}$ the intersection multiplicity of the curves $f=0$ and $g=0$ at the origin. Then $i_0(f,g)\geqslant(\ord_0 f)(\ord_0 g)$ with equality if and only if the curves $f=0$ and $g=0$ are transverse at $0$ i.e. don't have common tangent at $0$. The curve $f=0$ has an isolated singular point at $0$ if $\ord_0 f>1$ and $\mu_0(f)<+\infty$. Note that $f$ is a homogeneous polynomial of degree $d>0$ if and only if $\ord_0 f=d$

\section{Results}

We keep the notations introduced in Introduction. For any $a\in\rr$ we denote by $[a]$ the integral part of $a$.

The main result of this note is
\begin{Twierdzenie}\label{thm1.1}
Let $f=0$ be a curve of degree $d>1$ with an isolated singular point at $0\in\ccc^2$. Suppose that $\ord_0 f<d $. Then
$$
\mu_0(f)\leqslant(d-1)^2-\left[\frac{d}{2}\right].
$$
\end{Twierdzenie}
We prove Theorem \ref{thm1.1} in Section \ref{sec3}. The bound in the theorem is exact.

\begin{Example}
Let $d>1$ be an integer. 
$$
\mbox{Put}\quad f(x,y)=\left\{\begin{array}{l}
\displaystyle\phantom{x}\prod_{i=1}^{d/2}(x+ix^2+y^2)\quad \mbox{ if\ \ } d\equiv 0 \pmod 2,\\

\displaystyle x\prod_{i=1}^{\frac{d-1}{2}}(x+ix^2+y^2)\quad \mbox{ if\ \ } d\not\equiv 0 \pmod 2.
\end{array}
\right.
$$
Then $f$ is a polynomial of degree $d$ and $\mu_0(f)=(d-1)^2-\left[\frac{d}{2}\right]$.
\end{Example}

\begin{Uwaga}
S.~M.~Gusein-Zade and N.~N.~Nekhoreshev using topological methods proved in \cite{Gusein} (Proposition 1) that if $\ord_0 f=2$ then $\mu_0 (f)\leqslant(d-1)^2-\left[\frac{d}{2}\right]\left(\left[\frac{d}{2}\right]-1\right)$. Another bound for the Milnor number follows from the Abhyankar-Moh theory of approximate roots (see \cite{Garcia}, Corollary 6.5). Suppose that the curve $f=0$ is unibranch at $0$ (i.e. $f$ is irreducible in the ring of formal power series $\ccc[[x,y]]$) and the unique tangent to $f=0$ at $0$ intersects the curve with multiplicity $d$. Then $\mu_0(f)\leqslant(d-1)^2-\left(\frac{d}{d_1}-1\right)(d-\ord_0f)$, where $d_1=\gcd(\ord_0f,d)$.
\end{Uwaga}
\begin{Twierdzenie}\label{thm1.4}
Let $f$ be a polynomial of degree $d>2$, $d\neq 4$. Then the following two conditions are equivalent
\begin{itemize}
\item[(i)] The curve $f=0$ passes through the origin and $\mu_0(f)=(d-1)^2-\left[\frac{d}{2}\right]$,
\item[(ii)] The curve $f=0$ has $d-\left[\frac{d}{2}\right]$ irreducible components. Each irreducible component of the curve passes through the origin. If $d\equiv 0\pmod 2$ then all components are of degree $2$ and intersect pairwise at $0$ with multiplicity~$4$. If $d\not\equiv 0\pmod 2$ then all but one component are of degree $2$ and intersect pairwise at $0$ with multiplicity $4$, the remaining component is linear and is tangent to all components of degree $2$.
\end{itemize}
\end{Twierdzenie}
The proof of Theorem \ref{thm1.4} is given in Section \ref{sec4}.

\begin{Uwaga}
The implication $(ii)\Rightarrow(i)$ holds for any $d>2$. The assumption $d\neq 4$ is necessary for the implication $(i)\Rightarrow(ii)$. Take $f(x,y)=x(y^3-x^2)$. Then $f$ is of degree $d=4$, $\mu_0(f)=(d-1)^2-\left[\frac{d}{2}\right]=7$ and the condition $(ii)$ fails.
\end{Uwaga}

\section{Preparatory lemmas}

Let us begin with the following well-known properties of the Milnor number.

\begin{Lemat}\label{lem2.1}\mbox{}
\begin{itemize}
\item[(i)] If $f=f_0\tilde{f}$ in $\ccc[x,y]$ with $f_0(0)\neq 0$ then $\mu_0(f)=\mu_0(\tilde{f})$.
\item[(ii)] If $f=f_1\cdots f_m$ in $\ccc[x,y]$, $f_i(0)=0$ for $i=1,\dots,m$ then 
$$
\mu_0(f)+m-1=\sum_{i=1}^m \mu_0(f_i)+2\sum_{1\leqslant i<j\leqslant m}i_0(f_i,f_j).
$$
\end{itemize}
\end{Lemat}

\noindent\emph{Proof.} See, e.g. \cite{Cassou}, Property 5.4.\mbox{}\hfill$\Box$
\begin{Lemat}\label{lem2.2}
Let $f$ be an irreducible polynomial, $f(0)=0$, of degree $d>1$. Then $\mu_0(f)\leqslant (d-1)(d-2)$ with equality if and only if the curve $f=0$ is rational, its projective closure $C$ has exactly one singular point $0\in\ccc^2$ and $f=0$ is unibranch at $0$.
\end{Lemat}
\noindent\emph{Proof.} Apply the formula for the genus $g$ of $C$ (\cite{Wall}, Corollary 7.1.3): $2g=(d-1)(d-2)-\sum_P(\mu_P+r_P-1)$, where $r_P$ is the number of branches of $C$ passing through $P$.\mbox{}\hfill$\Box$

\begin{Example}
The polynomial $f=x^{d-1}+y^d$ is irreducible and $\mu_0(f)=(d-1)(d-2)$.
\end{Example}

Let $f=0$ be a curve of degree $d>1$ with an isolated singular point at $0\in\ccc^2$. Let $f_i=0$, $i=1,\dots,m$ be irreducible components of $f=0$ passing through $0$ and let $d_i=\deg f_i$ for $i=1,\dots,m$. Then $f=f_0f_1\cdots f_m$ in $\ccc[x,y]$, where $f_0(0)\neq 0$.

In what follows we keep the assumptions introduced above.

\begin{Lemat}\label{lem2.4}
Let $\Lambda=\{(i,j)\in\nn^2:1\leqslant i<j\leqslant m,\mbox{ the curves }f_i=0,\mbox{ } f_j=0 \mbox{ are tranverse and }d_i>1\mbox{ or }d_j>1\}$. Then $\mu_0(f)\leqslant (d-1)^2-d+m-2(\sharp\Lambda)$.
\end{Lemat}
\noindent\emph{Proof.} Let $\tilde{f}=f_1\cdots f_m$. Observe that for $(i,j)\in\Lambda$ we have $i_0(f_i,f_j)=(\ord_0f_i)(\ord_0f_j)<d_id_j$ since $d_i>1$ or $d_j>1$ (if $f_i$ is irreducible of degree $d_i> 1$ then $\ord_0f_i<d_i$). By Lemma \ref{lem2.2} we get $\mu_0(f_i)\leqslant(d_i-1)(d_i-2)$ for $i=1,\dots,m$. Now, Lemma \ref{lem2.1} implies
\begin{eqnarray*}
\lefteqn{\displaystyle\mu_0(f)+m-1=\mu_0(\tilde{f})+m-1\leqslant}\\
&&\displaystyle\leqslant\sum_{i=1}^m(d_i-1)(d_i-2)+
2\sum_{(i,j)\not\in\Lambda}i_0(f_i,f_j)+2\sum_{(i,j)\in\Lambda}i_0(f_i,f_j)\leqslant\\
&&\displaystyle\leqslant\sum_{i=1}^m(d_i-1)(d_i-2)+2\sum_{(i,j)\not\in\Lambda}d_id_j+
2\sum_{(i,j)\in\Lambda}(d_id_j-1)=\\
&&=(d-1)^2-d+2m-2(\sharp\Lambda)-1
\end{eqnarray*}
and the lemma follows.\mbox{}\hfill$\Box$

\begin{Lemat}\label{lem2.5}
We have that $\mu_0(f)\leqslant(d-1)^2-d+m$. The equality $\mu_0(f)=(d-1)^2-d+m$ holds if and only if $\mu_0(f_i)=(d_i-1)(d_i-2)$ and $i_0(f_i,f_j)=d_id_j$ for $1\leqslant i<j\leqslant m$.
\end{Lemat}
\noindent\emph{Proof.} The inequality $\mu_0(f)\leqslant(d-1)^2-d+m$ follows from Lemma \ref{lem2.4} (see also \cite{Gwozdziewicz}, Proposition 6.3). By Lemma \ref{lem2.1} we can rewrite the equality $\mu_0(f)=(d-1)^2-d+m$ in the form
$$
\sum_{i=1}^m\mu_0(f_i)+2\sum_{1\leqslant i<j\leqslant m}i_0(f_i,f_j)=\sum_{i=1}^m(d_i-1)(d_i-2)+2\sum_{1\leqslant i<j\leqslant m}d_id_j.
$$
We have that $\mu_0(f_i)\leqslant(d_i-1)(d_i-2)$ by Lemma \ref{lem2.2} and $i_0(f_i,f_j)\leqslant d_id_j$ by B\' ezout's theorem which together with the equality above imply the lemma.\mbox{}\hfill$\Box$ 

\begin{Lemat}\label{lem2.6}
Let $d>2$. If $\mu_0(f)=(d-1)^2-d+m$ then all irreducible components of the curve $f=0$ pass through $0\in\ccc^2$.
\end{Lemat}
\noindent\emph{Proof.} Let $\tilde{d}=\deg\tilde{f}$. Clearly $\tilde{d}\leqslant d$. We have $(d-1)^2-d+m=\mu_0(f)=\mu_0(\tilde{f})$ and $\mu_0(\tilde{f})\leqslant(\tilde{d}-1)^2-\tilde{d}+m$ by Lemma \ref{lem2.5}. The inequalities $d>2$, $\tilde{d}\leqslant d$ and $(d-1)^2-d\leqslant(\tilde{d}-1)^2-\tilde{d}$ imply $\tilde{d}=d$. Therefore $\deg f_0=d-\tilde{d}=0$ and $f_0$ is a constant.\mbox{}\hfill$\Box$

\begin{Lemat}\label{lem2.7}
$\quad\sharp\{i\in[1,m]:d_i>1\}\leqslant d-m$.
\end{Lemat}
\noindent\emph{Proof.} Let $I=\{i\in[1,m]:d_i>1\}$. Then $I=\{i\in[1,m]:\ord_0f_i<d_i\}$ and $\sharp I\leqslant \sum_{i\in I}(d_i-\ord_0 f_i)=\sum_{i=1}^m(d_i-\ord_0f_i)\leqslant d-\ord_0f\leqslant d-m$.\mbox{}\hfill$\Box$\vspace{2ex}

A line $l=0$ is tangent to the curve $f=0$ (at $0\in\ccc^2$) if $i_0(f,l)>\ord_0f$. We denote by $T(f)$ the set of all tangents at $0$ to $f=0$. We have $\sharp T(f)\leqslant\ord_0f\leqslant\deg f-1$ if $f$ is not homogeneous. For two polynomials $f,g$: $T(fg)=T(f)\cup T(g)$. Therefore we get $\sharp T(fg)\leqslant\sharp T(f)+\sharp T(g)-1$ if $T(f)\cap T(g)\neq\emptyset$.

\begin{Lemat}\label{lem2.8}
Let $f_i$, $i=1,\dots,k$ be irreducible polynomials of degree $d_i>1$ such that $f_i(0)=0$ for $i=1,\dots,k$. Suppose that the curves $f_i=0$, $i=1,\dots,k$ have a common tangent at $0$. Then $\sharp T(f_1\cdots f_k)\leqslant\sum_{i=1}^k(d_i-1)-k+1$.
\end{Lemat}
\noindent\emph{Proof.} If $k=1$ it is clear. Suppose that $k>1$ and that the lemma is true for the sequences of $k-1$ polynomials. Let $f_1,\dots,f_k$ be a sequence of $k$ irreducible polynomials of degree $>1$ such that the curves $f_i=0$, $i=1,\dots,k$ have a common tangent. Then by the induction hypothesis $\sharp T(f_1\cdots f_{k-1})\leqslant\sum_{i=1}^{k-1}(d_i-1)-(k-1)+1$. On the other hand $\sharp T(f_k)\leqslant d_k-1$ and we get $\sharp T(f_1\cdots f_k)\leqslant\sharp T(f_1\cdots f_{k-1})+\sharp T(f_k)-1\leqslant\sum_{i=1}^k(d_i-1)-(k-1)$, since $f_1\cdots f_{k-1}$ and $f_k$ have a common tangent.\mbox{ }\hfill$\Box$

\section{Proof of Theorem \ref{thm1.1}\label{sec3}}

Let $f$ be a polynomial of degree $d>1$ such that $f(0)=0$ and $\mu_0(f)<+\infty$. We assume that $f$ is not homogeneous. Let $m$ be the number of irreducible components of the curve $f=0$ passing through $0\in\ccc^2$.

\begin{Lemat}\label{lem3.1}
If $m\leqslant d-\left[\frac{d}{2}\right]$ then $\mu_0(f)\leqslant (d-1)^2-\left[\frac{d}{2}\right]$. The equality $\mu_0(f)=(d-1)^2-\left[\frac{d}{2}\right]$ implies $m=d-\left[\frac{d}{2}\right]$.
\end{Lemat}
\noindent\emph{Proof.} Suppose that $m\leqslant d-\left[\frac{d}{2}\right]$. By the first part of Lemma \ref{lem2.5} we get $\mu_0(f)\leqslant (d-1)^2-d+m\leqslant(d-1)^2-\left[\frac{d}{2}\right]$. If $\mu_0(f)=(d-1)^2-\left[\frac{d}{2}\right]$ then $(d-1)^2-d+m=(d-1)^2-\left[\frac{d}{2}\right]$, so $m=d-\left[\frac{d}{2}\right]$.\mbox{ }\hfill$\Box$

\begin{Lemat}\label{lem3.2}
\hspace{1em}If $m\geqslant d-\left[\frac{d}{2}\right]+1$ then $\mu_0(f)<(d-1)^2-\left[\frac{d}{2}\right]$.
\end{Lemat}
\noindent\emph{Proof.} Write $f=f_0f_1\cdots f_m$, where $f_0(0)\neq 0$ and $f_i$ are irreducible with $f_i(0)=0$ for $i=1,\dots,m$. Since the sequence $d\mapsto d-\left[\frac{d}{2}\right]$ is increasing it suffices to check the lemma for the polynomial $\tilde{f}=f_1\cdots f_m$. In what follows we write $f$ instead of $\tilde{f}$ and put $d_i=\deg f_i$ for $i=1,\dots,m$. Since $f$ is not homogeneous we have $m<d$. Let $k=\sharp\{i: d_i>1\}$. From $d_1+\cdots+d_m=d$ it follows that $k\geqslant 1$. Note also that $m-k=\sharp\{i:d_i=1\}>0$ since by Lemma \ref{lem2.7} we have $k\leqslant d-m$ and consequently $m-k\geqslant m-(d-m)=2m-d\geqslant d-2\left[\frac{d}{2}\right]+2$. We label $f_1,\dots,f_m$ so $d_1\geqslant\dots\geqslant d_m$. Therefore we get $d_1\geqslant\dots\geqslant d_k\geqslant 2$ and $d_{k+1}=\dots=d_m=1$. Let us consider two cases.\vspace{0.5em}

\noindent Case 1. The curves $f_1=0$, $\dots$, $f_k=0$ have a common tangent. By Lemma \ref{lem2.8} we have $\sharp T(f_1\cdots f_k)\leqslant\sum_{i=1}^k(d_i-1)-k+1=\sum_{i=1}^m(d_i-1)-k+1=d-m-k+1$. Therefore we get $m-k-\sharp T(f_1\cdots f_k)\geqslant m-k-(d-m-k+1)=2m-d-1$. Note that $2m-d-1\geqslant 2\left(d-\left[\frac{d}{2}\right]+1\right)-d+1>0$. Thus there are at least $2m-d-1>0$ linear forms in the sequence $f_{k+1},\dots,f_m$ that are transverse to the curve $f_1\cdots f_k=0$ and we get $\sharp\Lambda\geqslant k(2m-d-1\geqslant 2m-d-1)$. Consequently by Lemma \ref{lem2.4} we obtain $\mu_0(f)\leqslant(d-1)^2-d+m-2(2m-d-1)=(d-1)^2+d-3m+2\leqslant(d-1)^2+d-3\left(d-\left[\frac{d}{2}\right]+1\right)+2=(d-1)^2-\left[\frac{d}{2}\right]-2\left(d-\left[\frac{d}{2}\right]\right)-1<(d-1)^2-\left[\frac{d}{2}\right]$.
\vspace{0.5em}

\noindent Case 2. The curves $f_1=0$, $\dots$, $f_k=0$ have no common tangent. Then for every linear form $f_j$, $k+1\leqslant j\leqslant m$ there exists a polynomial $f_i$, $1\leqslant i\leqslant k$ such that $f_i$, $f_j$ are transverse. Therefore $\sharp\Lambda\geqslant m-k$ and by Lemma \ref{lem2.4} we get $\mu_0(f)\leqslant(d-1)^2-d+m-2(m-k)=(d-1)^2-d-m+2k$. Since by Lemma \ref{lem2.7} we have $k\leqslant d-m$ the above bound for $\mu_0(f)$ implies $\mu_0(f)\leqslant(d-1)^2+d-3m\leqslant(d-1)^2+d-3\left(d-\left[\frac{d}{2}\right]+1\right)=(d-1)^2-\left[\frac{d}{2}\right]-2\left(d-2\left[\frac{d}{2}\right]\right)-3<(d-1)^2-\left[\frac{d}{2}\right]$.

Now from Lemmas \ref{lem3.1} and \ref{lem3.2} we get $\mu_0(f)\leqslant(d-1)^2-\left[\frac{d}{2}\right]$ which proves Theorem \ref{thm1.1}. \mbox{ }\hfill$\Box$

\section{Proof of Theorem \ref{thm1.4}\label{sec4}}

\begin{Lemat}\label{lem4.1}
Let $f,g\in\ccc[x,y]$ be irreducible polynomials, $\deg f=3$, $\deg g=2$, $f(0)=g(0)=0$. Suppose that the curve $f=0$ has a singular point at $0$ and $\sharp T(f)=1$. Then
$$
i_0(f,g)<(\deg f)(\deg g)=6.
$$
\end{Lemat}
\noindent\emph{Proof.} If $f=0$ and $g=0$ have no common tangent then $i_0(f,g)=\ord_0f\ord_0g=\ord_0f<3$. Thus we may assume that $f=x^2+f^+$, $g=x+g^+$, where $f^+$, $g^+$ are homogeneous forms. We get $i_0(f,g)=i_0(f-xg,g)=i_0(f^+-xg^+,g)=3\cdot 1$ for $f^+-xg^+=0$ and $g=0$ have no common tangent.\mbox{ }\hfill$\Box$

\begin{Lemat}\label{lem4.2}
Let $f$ be a polynomial of degree $d>2$ such that $f(0)=0$. Suppose that $\mu_0(f)=(d-1)^2-d+m$, where $m$ is the number of irreducible components of the curve $f=0$ passing through $0\in\ccc^2$. Then $f=f_1\cdots f_m$ in $\ccc[x,y]$ with irreducible $f_i$, $f_i(0)=0$ for $i=1,\dots,m$. Let $d_i=\deg f_i$ for $i=1,\dots,m$. Then $\sharp T(f_i)=1$ and $i_0(f_i,f_j)=d_id_j$ for $i<j$. If $m<d$ then $f=0$ has at most one linear component and has no two components of degree $2$ and $3$.
\end{Lemat}
\noindent\emph{Proof.} By Lemmas \ref{lem2.5} and \ref{lem2.6} we get $f=f_1\cdots f_m$ in $\ccc[x,y]$, $f_i$ irreducible, $f_i(0)=0$, $i_0(f_i,f_j)=d_id_j$ and $\mu_0(f_i)=(d_i-1)(d_i-2)$. We have $\sharp T(f_i)=1$ for the curve $f_i=0$ has only one branch at $0$ by Lemma \ref{lem2.2}.

Suppose that the curve $f=0$ has two linear components $f_j=0$ and $f_k=0$, $j\neq k$. Then there is no component $f_i=0$ of $f=0$ of degree $d_i>1$ (if $f_i=0$ had degree $d_i>1$ then we would get $i_0(f_i,f_j)=i_0(f_i,f_k)=d_i>1$ which is impossible for $\sharp T(f_i)=1$). Therefore if $f=0$ has two linear components then all components are linear and intersect pairwise with multiplicity $1$. Thus $m=d$ and $f$ is a homogeneous form of degree $d$.

Therefore, if $m<d$ then there exists at most one linear component. Since $i_0(f_i,f_j)=d_id_j$ if $i<j$ there are no two components of degree $2$ and $3$ by Lemma \ref{lem4.1}.\mbox{ }\hfill$\Box$\vspace{1em}

Now, we can pass to the proof of Theorem \ref{thm1.4}\vspace{0.5em}

\noindent$(i)\Rightarrow(ii)$\\
Assume that $\mu_0(f)=(d-1)^2-\left[\frac{d}{2}\right]$, where $d>1$ and $d\neq 4$. Then $f=0$ has $m=d-\left[\frac{d}{2}\right]$ irreducible components passing through $0\in\ccc^2$ by Lemmas \ref{lem3.1} and \ref{lem3.2}. We have $\mu_0(f)=(d-1)^2-\left[\frac{d}{2}\right]=(d-1)^2-d+m$ with $m=d-\left[\frac{d}{2}\right]<d$ and by Lemma \ref{lem4.2} we can write $f=f_1\cdots f_m$, $f_i\in\ccc[x,y]$ irreducible, $f_i(0)=0$. We label $f_1,\dots,f_m$ so that $d_1\geqslant\dots\geqslant d_m\geqslant 1$.
\vspace{0.5em}

\noindent Case 1. $\quad d\equiv 0\pmod 2$.\\
Then $m=\frac{d}{2}$ and $d_1+\cdots+d_m=d$. This is possible if and only if $(d_1,\dots,d_m)=(2,\dots,2)$ or $(d_1,\dots,d_m)=(3,2,\dots,2,1)$, where $2$ appears $m-2=\frac{d}{2}-2>0$ times since $d>4$. If $(d_1,\dots,d_m)=(2,\dots,2)$ then the theorem follows from Lemma \ref{lem4.2}. The case $(d_1,\dots,d_m)=(3,2,\dots,2,1)$ can not occur by Lemma \ref{lem4.1}
\vspace{0.5em}

\noindent Case 2. $\quad d\not\equiv 0\pmod 2$.\\
In this case we have $m=\frac{d+1}{2}$. From $d_1+\cdots+d_m=d$ it follows that $(d_1,\dots,d_m)=(2,\dots,2,1)$. We apply Lemma \ref{lem4.2}. \vspace{0.5em}

The implication $(ii)\Rightarrow(i)$ follows immediately from Lemma \ref{lem2.1} $(ii)$.\mbox{ }\hfill$\Box$

\vspace{1ex}
\noindent Department of Mathematics\\
Kielce University of Technology\\
Al. 1000 L PP 7\\
25-314 Kielce, Poland\\
e-mail: matap@tu.kielce.pl
\end{document}